\documentstyle[a4wide]{article}
\begin{document}
\newcommand{\qed}{\nopagebreak
\begin{flushright} {\em q.e.d.} \end{flushright}}
\newcommand{\br}{\mbox{$\overrightarrow{\mbox{\bf b}}$}}
\newcommand{\bl}{\mbox{$\overleftarrow{\mbox{\bf b}}$}}
\newcommand{\kr}{\mbox{$\overrightarrow{\mbox{\bf k}}$}}
\newcommand{\kl}{\mbox{$\overleftarrow{\mbox{\bf k}}$}}
\newcommand{\ms}{\mbox{\boldmath{${\sigma}$}}}
\newcommand{\md}{\mbox{\boldmath{${\delta}$}}}
\newcommand{\mal}{\mbox{\boldmath{${\alpha}$}}}
\newcommand{\mbe}{\mbox{\boldmath{${\beta}$}}}
\newcommand{\mep}{\mbox{\boldmath{${\epsilon}$}}}
\newcommand{\met}{\mbox{\boldmath{${\eta}$}}}
\newcommand{\mtv}{\mbox{\boldmath{$T$}}}
\newcommand{\tv}{\mbox{\rm T}}
\newcommand{\msv}{\mbox{\boldmath{$S$}}}
\newcommand{\mrv}{\mbox{\boldmath{$R$}}}
\newcommand{\mj}{\mbox{\bf 1}}
\newcommand{\ccc}{\mbox{\bf CartCl}}
\newcommand{\frd}{\mbox{\bf Finord}}
\newcommand{\cc}{\mbox{\bf Cart}}
\newcommand{\sym}{\mbox{\bf SyMon}}
\newcommand{\symcl}{\mbox{\bf SyMonCl}}
\newcommand{\skup}{\mbox{\bf Set}}
\newcommand{\mc}{\mbox{\bf c}}
\newcommand{\mw}{\mbox{\bf w}}
\newcommand{\mk}{\mbox{\bf k}}
\newcommand{\mb}{\mbox{\bf b}}
\newcommand{\ri}{\mbox{\rm I}}
\newcommand{\mi}{\mbox{\scriptsize{\rm I}}}
\newcommand{\mpm}{\mbox{\boldmath{$\pi$}}}           
\newcommand{\dkz}{{\bf Proof:}{\hspace{0.3cm}} \nopagebreak}
\newcommand{\str}{\rightarrow}

\baselineskip=1.05\baselineskip

\title 
{The Maximality of Cartesian Categories}

\author{ }

\date{ }
\maketitle
\begin{abstract}
\noindent It is proved that equations between arrows assumed for cartesian
categories are maximal in the sense that extending them with any new equation
in the language of free cartesian categories collapses a cartesian category
into a preorder. An analogous result holds for categories with binary
products, which may lack a terminal object. The proof is based on a coherence
result for cartesian categories, which is related to model-theoretical
methods of normalization.
\end{abstract}
\vspace{0.5cm}
\noindent The equations between arrows assumed for cartesian categories
are maximal in the sense that extending them with new equations
collapses the categories into preorders (i.e. categories in which between
any two objects there is at most one arrow). The equations envisaged for
the extension are in the language of free cartesian categories generated 
by sets of objects, and variables for arrows don't occur in them. If 
such an equation doesn't hold in the free cartesian category generated by 
a set of objects, then any cartesian category in which this equation holds
is a preorder. An
analogous result is provable for categories with binary products, which
differ from cartesian categories in not necessarily having a terminal
object.

The proof of these results, which we are going to present below, is based on a
coherence property of cartesian categories. This coherence, which is
ultimately inspired by the geometric modelling of categories of \cite{GFC},
is related to
model-theoretic methods of normalization. It permits to establish uniqueness
of normal form for arrow terms without proceeding via the Church-Rosser
property for reductions. It also yields an easy decision procedure
for the commuting of diagrams in free cartesian categories.

\section{Cartesian categories}
A {\em category with binary products} is a category with a binary operation
$\times$  on objects, {\em projection arrows}
\[ \mk^{1}_{A,B}:A \times B \str A, \]
\[ \mk^{2}_{A,B}:A \times B \str B, \]
and the {\em pairing operation} that to a pair of arrows ($f:C \str A$,
$g:C \str B$) assigns the arrow $\langle f,g\rangle:C \str A \times B$.
The arrows must satisfy the equations
\begin{tabbing}
\hspace{15em} \= $(\beta 1)$ \hspace{2em} \= $\mk^{1}_{A,B} \circ
\langle f,g \rangle = f$, \\[0.3cm]
\> $(\beta 2)$ \> $\mk^{2}_{A,B} \circ
\langle f,g \rangle = g$, \\[0.3cm]
\> $(distr)$ \> $\langle f,g \rangle \circ h =
\langle f \circ h , g \circ h \rangle$, \\[0.3cm]
\> $(\eta)$ \> $\langle \mk^{1}_{A,B} , \mk^{2}_{A,B} \rangle = \mj_{A
\times B}$.
\end{tabbing}

A category has a {\em terminal object} T iff it has the special arrows
\[\mk_{A} : A \str \tv , \]
which satisfy the equation
\[ (\mk){\mbox{\hspace{2em}}}{\mbox{\rm for }}f:A \str \tv,
{\mbox{\hspace{2em}}}f=\mk_{A}.\]

A {\em cartesian category} is a category with binary products and a
terminal object.

In terms of the projection arrows and the pairing operation on arrows we may
define in every category with binary products the following arrows:
\begin{tabbing}
\hspace{5em} \= $\br_{A,B,C}$ \= $=_{def.}$ \= ${\langle}{\langle}\mk^{1}_{A,B \times C},
\mk^{1}_{B,C} \circ \mk^{2}_{A,B \times C}{\rangle},\mk^{2}_{B,C} \circ
\mk^{2}_{A,B \times C}{\rangle} $ \\[0.1cm]
\> \> \> of type $A \times (B \times C) \str (A \times B) \times C$, \\[0.2cm]
\> $\bl_{A,B,C}$ \> $=_{def.}$ \> ${\langle}\mk^{1}_{A,B} \circ \mk^{1}_{A \times B,C},
{\langle}\mk^{2}_{A,B} \circ \mk^{1}_{A \times B,C},\mk^{2}_{A \times B,C}
{\rangle}{\rangle}$ \\[0.1cm]
\> \> \> of type $(A \times B) \times C \str A \times (B \times C)$, \\[0.2cm]
\> $\mc_{A,B}$ \> $=_{def.}$ \> ${\langle}\mk^{2}_{A,B},\mk^{1}_{A,B}{\rangle}$
\\[0.1cm]
\> \> \> of type $A \times B \str B \times A$, \\[0.2cm]
\> $\mw_{A}$ \> $=_{def.}$ \> ${\langle}\mj_{A},\mj_{A}{\rangle}$ \\[0.1cm]
\> \> \> of type $A \str A \times A$.
\end{tabbing}

We may also define the {\em product operation on arrows}, which to a
pair of arrows ($f:A \str B$, $g:C \str D$) assigns the arrow
$f \times g:A \times C \str B \times D$:
\[ f \times g =_{def.} \langle f \circ \mk^{1}_{A,C} , g \circ \mk^{2}_{A,C}
\rangle .\]
In every cartesian category we also have the arrows
\begin{tabbing}
\hspace{5em} \= $\br_{A,B,C}$ \= $=_{def.}$
\= ${\langle}{\langle}\mk^{1}_{A,B \times C},
\mk^{1}_{B,C} \circ \mk^{2}_{A,B \times C}{\rangle},\mk^{2}_{B,C} \circ
\mk^{2}_{A,B \times C}{\rangle} $ \kill
\> $\ms_{A}$ \> $=_{def.}$ \> $\langle \mk_{A} , \mj_{A} \rangle$ \\[0.1cm]
\> \> \> of type $A \str \tv \times A$, \\[0.2cm]
\> $\md_{A}$ \> $=_{def.}$ \> $\langle \mj_{A} , \mk_{A} \rangle$ \\[0.1cm]
\> \> \> of type $A \str A \times \tv$.
\end{tabbing}

These definitions of category with binary products and cartesian category
are equivalent to standard definitions (see \cite{Lambek}, Part 0, Chapter 5,
and Part I, Chapter 3), where
\[A\stackrel{\mk^1_{A,B}}{\longleftarrow} A\times B \stackrel{\mk^2_{A,B}}
{\longrightarrow} B \]
is a product diagram with the desired universal property.

\section{Graphs of arrow terms in free cartesian categories}
Consider the free cartesian category \cc \ generated by a set
of objects $\cal O$ called {\em letters} (T is not a letter).
This category is the image of $\cal O$ under the left adjoint to
the forgetful functor from the category of
cartesian categories, with cartesian structure-preserving
functors as arrows, to the category
of sets (of objects). The construction of \cc \ out of
syntactic material is explained in detail in \cite{Lambek} (Part I,
Chapter 4; note that there the name ``\cc'' has a different meaning).

For
an object $A$ of \cc \ let the {\em letter length} $|A|$ of $A$ be the number
of occurrences of letters in $A$. For example, if $p$ and $q$
are letters, then
$|((p \times q) \times p) \times (\tv \times p)|$ is 4.
Let $f:A \str B$ be an arrow term
of \cc \ and let $|A|=n$ and $|B|=m$. To $f$ we associate a
function $\Gamma_{f}$ from $\{1, \ldots ,m\}$ to $\{1, \ldots ,n\}$,
called the {\em graph} of $f$. If $m=0$, then $\{1, \ldots ,m\}$ is $\emptyset$. The function $\Gamma_{f}$ is defined by induction
on the complexity of $f$ in the following manner.

If $f$ is of the form $\mj_{A}:A \str A$, then $m=n$ and $\Gamma_{f}(i)=i$,
where $i \in \{1, \ldots ,m\}$. If $f$ is of the form $\mk^{1}_{A,B}:
A \times B \str A$, then $\Gamma_{f}(i)=i$, and if $f$ is of the form
$\mk^{2}_{A,B}:A \times B \str B$, then $\Gamma_{f}(i)=i+|A|$. If $f$ is
of the form $\mk_{A}:A \str \tv$, then $\Gamma_{f}$ is the empty function
(i.e., $\Gamma_{f}$ is $\emptyset: \emptyset \str \{1, \ldots ,|A|\}$).

If $f$ is of the form $\langle g,h \rangle: C \str A \times B$, with
$g:C \str A $ and $h:C \str B$, then
for $i \leq |A|$ we have $\Gamma_{f}(i)=\Gamma_{g}(i)$
and
for $i > |A|$ we have $\Gamma_{f}(i)=\Gamma_{h}(i)$. Finally, if 
$f$ is of the form $h \circ g$, then $\Gamma_{f}(i)=\Gamma_{g}(\Gamma_{h}
(i))$.

For $f: A \str B$, the graph $\Gamma_{f}$ can be interpreted as connecting
an occurrence of a letter $p$ in $A$ with a finite set of occurrences of
$p$ in $B$ (this set may have more than one member, it may be a singleton,
or it may be empty).

We can establish the following lemma.
\\[0.3cm]
{\bf Lemma 1}\hspace{1em}
{\em If $f=g$ in \cc, then $\Gamma_{f}=\Gamma_{g}$.}
\\[0.3cm]
\dkz
We proceed by induction on the length of the derivation of $f=g$ in
\cc. In this induction it is essential to check that
for the equations ($\beta 1$), ($\beta 2$), ($distr$), ($\eta$) and
(\mk), the graphs of the two sides of the equation must be equal. The
induction step, which involves the rules of symmetry and transitivity
of equality, as well as congruence with composition and pairing, is
quite trivial.
\qed
We shall demonstrate the converse of this lemma in Lemma 3 below. 
\section{Normal form of arrow terms in free cartesian categories}
From now on, identity of arrow terms of \cc \ will be taken up
to associativity of composition. So, for example,
$\mk^{1}_{B,C} \circ (\mk^{2}_{A,B \times C} \circ
\mj_{A \times (B \times C)})$
will be considered to be the {\em same} arrow term as
$(\mk^{1}_{B,C} \circ \mk^{2}_{A,B \times C}) \circ
\mj_{A \times (B \times C)}$,
and we may omit parentheses in compositions. (Formally, we may work
with equivalence classes of arrow terms.)

An arrow term of \cc \ is called an {\em atomized} \mk{\em -composition}
iff it is of the form $f_{n} \circ \cdots \circ f_{1} : A \str B$, with
$n \geq 1$, where $B$ is a letter and each $f_{i}$ is either of the form
$\mk^{1}_{C,D}$ or of the form $\mk^{2}_{C,D}$.

Arrow terms of \cc \ in {\em normal form} are defined inductively as follows:
\begin{enumerate}
\item every arrow term $\mj_{A}$ with $A$ a letter is in normal form;
\item every atomized \mk-composition is in normal form;
\item every arrow term $\mk_{A}$ is in normal form;
\item if $f:C \str A$ and $g:C \str B$ are in normal form,
then $\langle f,g \rangle$ is
in normal form.
\end{enumerate}
The arrow terms defining \br, \bl, \mc, \mw, \ms\ and \md\
arrows in section 1
are in normal form if $A$, $B$ and $C$ are letters.

For the normal form of arrow terms of \cc \ and their graphs we can prove the
following fundamental lemma.
\\[0.3cm]
{\bf Lemma 2}\hspace{1em}
{\em Suppose $f,g:A \str B$ are arrow terms of \cc \ in normal form.
Then $\Gamma_{f}=\Gamma_{g}$ iff $f$ and $g$ are the same arrow term.}
\\[0.3cm]
\dkz
Suppose $f,g:A \str B$ are different arrow terms of \cc \ in normal form.
We shall show that in that case $\Gamma_{f} \not= \Gamma_{g}$ by
induction on the number of pairing operations in $f$.

If this number is zero, then it must be zero too in $g$, because $B$
is a letter or \tv. If $f$ is of the form $\mj_{A}$, then $g$ must also
be $\mj_{A}$, and if $f$ is of the form $\mk_{A}$, then $g$ must also be
$\mk_{A}$. So $f$ cannot be $\mj_{A}$ or $\mk_{A}$. The only
remaining possibility is that $f$ be an atomized \mk-composition of
the form $f_{n} \circ \cdots \circ f_{1}$. Then $g$ must be an atomized
\mk-compsition, too; let it be of the form $g_{m} \circ \cdots \circ
g_{1}$. It is excluded that $f$ be of the form
\[f_{n} \circ \cdots \circ f_{n-k} \circ g_{m} \circ \cdots \circ g_{1}\]
for $n-k \geq 1$, because the codomain of $g_{m}$ is a letter.
Analogously, it is excluded that $g$ be of the form
\[g_{m} \circ \cdots \circ g_{m-k} \circ f_{n} \circ \cdots \circ f_{1}\]
for $m-k \geq 1$.

So for some $i$ we must have that $f_{i}$ and $g_{i}$
are different; let $j$ be the least such $i$. Then one of $f_{j}$ and
$g_{j}$ is $\mk^{1}_{C,D}$ while the other is $\mk^{2}_{C,D}$, and the
letter $B$ must occur in both $C$ and $D$. It follows easily that
$\Gamma_{f} \not= \Gamma_{g}$.

If the number of pairing operations in $f$ is at least one, then $f$
is of the form $\langle f_{1},f_{2} \rangle: A \str B_{1} \times B_{2}$,
and hence $g$, which is in normal form, must be of the form 
$\langle g_{1},g_{2} \rangle$. Since $f$ is different from $g$, either
$f_{1}$ and $g_{1}$ or $f_{2}$ and $g_{2}$ are two different arrow terms
with identical domains and codomains, which are both in normal form.
By the induction hypothesis,
either $\Gamma_{f_{1}} \not= \Gamma_{g_{1}}$ or
$\Gamma_{f_{2}} \not= \Gamma_{g_{2}}$, and in both cases we can infer
$\Gamma_{f} \not= \Gamma_{g}$. This concludes the induction.

So, by contraposition, if $\Gamma_{f}=\Gamma_{g}$, then $f$ and $g$
are the same arrow term. Since the converse implication is trivial,
this proves the lemma.
\qed

Every arrow term of \cc \ can be reduced to normal form by using the
following reductions:

\begin{center}
\rm
\begin{picture}(300,240)
\put(80,20){\makebox(0,0){$g:C \str \tv$}}
\put(220,20){\makebox(0,0){$\mk_{C}$}}
\put(80,40){\makebox(0,0){$f:C \str A \times B$}}
\put(220,40){\makebox(0,0){$\langle \mk^{1}_{A,B} \circ f , \mk^{2}_{A,B}
\circ f \rangle$}}
\put(150,60){\makebox(0,0){\em atomizing reductions}}
\put(80,80){\makebox(0,0){$\langle f , g \rangle \circ h$}}
\put(220,80){\makebox(0,0){$\langle f \circ h , g \circ h \rangle$}}
\put(80,100){\makebox(0,0){$\mk^{2}_{A,B} \circ \langle f , g \rangle$}}
\put(220,100){\makebox(0,0){$g$}}
\put(80,120){\makebox(0,0){$\mk^{1}_{A,B} \circ \langle f , g \rangle$}}
\put(220,120){\makebox(0,0){$f$}}
\put(150,140){\makebox(0,0){\em pairing reductions}}
\put(80,160){\makebox(0,0){$\mj_{B} \circ f$}}
\put(220,160){\makebox(0,0){$f$}}
\put(80,180){\makebox(0,0){$f \circ \mj_{A}$}}
\put(220,180){\makebox(0,0){$f$}}
\put(150,200){\makebox(0,0){\em identity reductions}}
\put(80,220){\makebox(0,0){redexes}}
\put(220,220){\makebox(0,0){contracta}}
\end{picture}
\end{center}

We reduce an arrow term $t_0$ of \cc\ to normal form via a sequence $t_0, t_1,
\ldots, t_n$, where $t_{i+1}$ is obtained from $t_i$ by a reduction step,
which is a replacement of a subterm of $t_i$ that is a redex by the
corresponding contractum. We call each $t_i$, $i<n$, in this sequence
a {\em candidate for reduction}.

To obtain that {\em every} sequence of
reduction steps terminates we must
exclude in the first kind of atomizing reduction that the
redex $f$ be of the form $\langle f_{1},f_{2} \rangle$.
We should exclude as well that it reduces to this form by other
reductions. We must also ensure that the respective occurrence
of $f$ does not
belong to a subterm of the form $h \circ f$ where $h$ is either
$\mk^{1}_{A,B}$ or $\mk^{2}_{A,B}$, or $h$ reduces to
$\mk^{1}_{A,B}$ or $\mk^{2}_{A,B}$. In the second kind of
atomizing reduction the redex $g$ should be different from $\mk_{C}$.
All this will be guaranteed if we add to the atomizing reductions
the provisos we are going to formulate below.

For an arrow term $t$ of \cc\, let $\gamma(t)$ be defined inductively
as follows:
\[\begin{array}{l}
\gamma(\mk_A)=2, \\
\gamma(t)=3 {\mbox{\rm\hspace{1em}if $t$ is $\mj_A$, $\mk^1_{A,B}$ or
$\mk^2_{A,B}$, }} \\
\gamma(g \circ f)= \gamma(g) \cdot \gamma(f), \\
\gamma(\langle f, g \rangle)= \gamma(f) + \gamma(g) + 1.
\end{array} \]

Let $\alpha(t)$ be $\gamma(t) \cdot n$, where $n$ is 1 if in $t$ there are
pairing operations $\langle \ldots, \ldots \rangle$ within the scope
of composition $\circ$; otherwise $n$ is 0.

If $h_2 \circ h_1$ is a subterm of an arrow
term $t$, then $h_1$ and $h_2$ are
called {\em compositional subterms} of $t$. Let a subterm $f$ of $t$ be called
{\em product-eliminative} (for this terminology, suggested by
natural deduction, see section 6) iff the pairing operation does not
occur in $f$, and $f$ is not a compositional subterm of $t$.

The proviso for the first atomizing reduction says that the
redex $f$ must be product-eliminative with respect to the candidate
for reduction $t$,
of which $f$ is a subterm, and $\alpha(t)=0$. The proviso for the second
atomizing reduction says just that $g$ is not $\mk_C$.

With these provisos, it is easy to
check that $t$ is in normal form iff there are no redexes
in it. (For that we rely on the fact that in \cc\ there is no arrow of type
$\tv \str A$ with $A$ a letter.)

For an object $A$ of \cc\ let the {\em length} of $A$ be the letter length
$|A|$ of $A$ plus the number of occurrences of $\times$ and T in $A$ (i.e.,
the length of $A$ is the number of occurrences of symbols in $A$).
For an arrow term $t$ of \cc\ let $\beta(t)$ be the sum of all the lengths of all
the targets of product-eliminative subterms of $t$. (If the same object is
the target of $n$ product-eliminative
subterms of $t$, then it is counted $n$ times in the
sum $\beta(t)$.)

Let the {\em degree} of a candidate for reduction
$t$ be the ordinal $\omega^2 \cdot \alpha(t) +
\omega \cdot \beta(t) + \gamma(t)$. It is easy to check
that by replacing a redex
of $t$ by a contractum the degree of the resulting arrow term
strictly decreases.
Then an induction up to $\omega^3$ (which is, of course, reducible to
an ordinary induction up to $\omega$) shows that every sequence of
reduction steps terminates.

Our provisos don't only
yield that every sequence of reduction steps terminates in a normal form,
but they optimize reductions in other
respects, too. For our purposes, however, it is enough to know that {\em some} sequence of reduction steps terminates in a normal form, so that the provisos for atomizing reductions are not essential. But the provisos do help to make it clear that such a sequence exists. All the reductions above (without the provisos, and hence with the provisos, too) are covered by equations of \cc.

Lemmata 1 and 2 guarantee that if $f=g$ is satisfied in \cc \ and $f'$ and
$g'$ are normal forms of $f$ and $g$, respectively, then $f'$ and $g'$
are the same arrow term. To prove that we conclude by Lemma 1
from $f=g$, $f=f'$ and $g=g'$, which we have in \cc, that $\Gamma_{f'}=
\Gamma_{g'}$. Then by Lemma~2, it follows that $f'$ and $g'$ are the
same arrow term.

So arrows of \cc \ have a unique normal form. Note that we have
demonstrated this without appealing to the Church-Rosser property
for our reductions.

\section{Coherence}

Before proving our theorem about the maximality of cartesian categories,
we establish a lemma converse to Lemma 1.
\\[0.3cm]
{\bf Lemma 3}\hspace{1em}
{\em Suppose $f,g:A \str B$ are arrow terms of \cc. If 
$\Gamma_{f}=\Gamma_{g}$,
then $f=g$ in \cc.}
\\[0.3cm]
\dkz
Suppose $f,g:A \str B$ are arrow terms of \cc \ such that
$\Gamma_{f}=\Gamma_{g}$. Then let $f'$ and $g'$ be the normal forms of
$f$ and $g$, respectively. Since $f=f'$ and $g=g'$ in \cc, by Lemma 1
we have $\Gamma_{f'}=\Gamma_{g'}$, and hence, by Lemma 2, $f'$ and $g'$
are the same arrow term. Then by the symmetry and transitivity of
equality it follows
that $f=g$ in \cc.
\qed

A notion of graph analogous to ours may be found in \cite{K1} (p. 94). 
A coherence result analogous to Lemma 3 is envisaged in \cite{K2} (p. 129),
and is demonstrated in \cite{MI} (Theorem 2.2), \cite{TS}
(Theorem 8.2.3, p. 207) and \cite{P}.

Lemmata 1 and 3 guarantee that there is a faithful cartesian
functor from \cc \  to the category $\frd^{op}$, whose objects are finite
ordinals and whose arrows are arbitrary functions from finite
ordinals to finite ordinals, with domains being targets and
codomains sources. This functor is onto on
objects and on arrows if \cc \ is generated by a nonempty set of letters.
The product of $\frd^{op}$, to which the product of \cc\ is
mapped, is simply addition, and the terminal object is zero.
If \cc\ is generated by a
single object, then $\frd^{op}$ is equivalent (but not isomorphic)
to \cc: it is the skeleton of \cc.

Our demonstration of uniqueness of the normal form in the preceding section,
which did not appeal to the Church-Rosser property of reductions,
is akin to model-theoretical methods of normalization
(see \cite{CUDS} and references therein). In the spirit of these methods, a
computation without reductions of the normal form of an arrow term $f$ of
\cc\ consists in finding $\Gamma_{f}$, and then constructing out of
$\Gamma_{f}$ an arrow term of \cc\ in normal form whose graph is
$\Gamma_{f}$.

The equivalence between $f=g$ and $\Gamma_{f}=\Gamma_{g}$, which
follows from Lemmata 1 and 3, provides an easy decision procedure
for commuting of diagrams in \cc. To check whether $f=g$ in \cc, we
just have to check whether $\Gamma_{f}=\Gamma_{g}$.
(An alternative decision procedure for the commuting of
diagrams in \cc\ can be obtained via reduction to normal form, according to
section 3.)

\section{The maximality theorem}
For arrow terms $f,g:A \str B$ of \cc, we say that $f=g$ {\em holds}
in a cartesian category $\cal C$ iff for every cartesian functor $F$
from \cc\ to $\cal C$ we have $F(f)=F(g)$ in $\cal C$.

Then we can prove our theorem.
\\[0.3cm]
{\bf Theorem}\hspace{1em}
{\em Suppose $f,g:A \str B$ are arrow terms of \cc\ such that in \cc\ we don't
have $f=g$. If $f=g$ holds in a cartesian category $\cal C$, then
$\cal C$ is a preorder.}
\\[0.3cm]
\dkz
Let $f,g:A \str B$ be arrow terms of \cc\ such that in \cc\ we don't have
$f=g$. By Lemma 3, it follows that $\Gamma_{f} \not= \Gamma_{g}$.
So there must be an occurrence of a letter $p$ in $B$ such that
if $p$ is the $i$-th letter symbol of $B$, counting from the left,
then $\Gamma_{f}(i) \not= \Gamma_{g}(i)$. Consider the substitution
instances $f',g':A' \str B'$ of $f$ and $g$ obtained by replacing every
letter by $p$. So the only letter in $A'$ and $B'$ is $p$. It is
clear that $\Gamma_{f'}(i) \not= \Gamma_{g'}(i)$.

Then there is an arrow term $h:p \times p \str A'$ of \cc\ made of
possibly multiple occurrences of the arrow terms
\mj, \br, \bl, \mc, \mw, \ms\ and \md, together with the product operation
on arrows and composition, such that $\Gamma_{h}(\Gamma_{f'}(i))=1$
and $\Gamma_{h}(\Gamma_{g'}(i))=2$.

There is also an arrow term $j:B' \str p$ that is either $\mj_{p}$
or an atomized \mk-composition (see section 3) such that
$\Gamma_{j}(1)=i$. Then $j \circ f' \circ h$ and $j \circ g' \circ h$
are both of type $p \times p \str p$, while
$\Gamma_{j \circ f' \circ h} = \Gamma_{\mk^{1}_{p,p}}$
and
$\Gamma_{j \circ g' \circ h} = \Gamma_{\mk^{2}_{p,p}}$.
Therefore, by Lemma 3, in \cc\ we have
$j \circ f' \circ h = \mk^{1}_{p,p}$ and
$j \circ g' \circ h = \mk^{2}_{p,p}$. So these two equations hold in
every cartesian category $\cal C$. If $f=g$ holds in $\cal C$,
then $f'=g'$ holds too, and so $\mk^{1}_{p,p} = \mk^{2}_{p,p}$
holds in $\cal C$.

Now suppose $s,t:C \str D$ are two arrows of $\cal C$. Then in $\cal C$
we have
\[ \mk^{1}_{D,D} \circ \langle s , t \rangle = \mk^{2}_{D,D} \circ 
\langle s , t \rangle,
\]
and so $s=t$, by ($\beta 1$) and ($\beta 2$).
\qed

An analogous maximality theorem can be proved for categories with
binary products. For that we have to replace \cc\ by the free category
with binary products $\cc^{-}$ generated by a set of letters. Graphs for
the arrow terms of $\cc^{-}$  are defined by just omitting the clause for
$\Gamma_{\mk_{A}}$ in section 2. In the definition of the normal
form in section 3 we omit clause (3), and the second kind of
atomizing reduction
in the same section is now superfluous. Exact analogues of
Lemmata 1, 2 and 3 are provable as before, while in the proof of the
Theorem we need not mention \ms\ and \md\ arrow terms in the second
paragraph. This way we show that, if $f,g:A \str B$ are arrow terms of
$\cc^{-}$ such that in $\cc^{-}$ we don't have $f=g$, and $f=g$
holds in a category with binary products $\cal C$, then $\cal C$ is
a preorder.

\section{A logical conclusion}

It is well known that equations between arrows in cartesian
categories correspond to equivalence between deductions in conjunctive
logic, including the constant true proposition. This equivalence
induces an equality relation on equivalence classes. As far as $\times$
is concerned, this equality between deductions permits us to
reduce every deduction into one in atomized normal form, where
elimination rules, which correspond to $\mk^{1}$ and $\mk^{2}$
arrows, precede introduction rules, which correspond to the
pairing operation, and the middle part between eliminations and
introductions is atomic. This atomized normal form corresponds to the
normal form of arrow terms of \cc\ in section 3.

The import of the maximality of cartesian categories for logic is
that the usual assumptions for equivalence between deductions in conjunctive
logic are optimal. These assumptions are wanted, because they are induced by
normalization of deductions, and no assumption is missing, because
any further equation would equate all deductions that share premises and
conclusions.

Of course, by duality, this solves the problem for disjunctive logic,
including the constant absurd proposition or not. This logic corresponds
to categories with binary coproducts, which when they have an initial
object are ``cocartesian'' categories, i.e. cartesian categories with 
arrows reversed.

It is natural to enquire whether similar maximality results hold for
other sorts of categories of interest to logic.  
An exactly analogous result for cartesian closed categories is proved 
in \cite{S1} (Theorem 1) and \cite{Maxccc}. However,
this result is independent
from the maximality result for cartesian categories proved in this 
note: neither result can be inferred from the other.

\end{document}